\documentclass[11pt]{amsart}

\usepackage{amsmath,amssymb,amsthm}
\usepackage[a4paper]{geometry}
\usepackage{graphicx}
\usepackage{color}
\usepackage{url}

\newtheorem{theorem}{Theorem}[section]
\newtheorem{lemma}[theorem]{Lemma}
\newtheorem{prop}[theorem]{Proposition}
\newtheorem{cor}[theorem]{Corollary}
\newtheorem{defi}{Definition}[section]
\newtheorem{rem}{Remark}[section]
\newtheorem{exple}{Example}[section]

\def\dlim{\stackrel{ (fdd)}{=}}
\newcommand{\rd}{{\mathbb R^d}}
\newcommand{\EM}{\ensuremath}

\newcommand{\dC}{\EM{\mathbb{C}}}
\newcommand{\cD}{\EM{\mathcal{D}}}
\newcommand{\dE}{\EM{\mathbb{E}}}
\newcommand{\dN}{\EM{\mathbb{N}}}
\newcommand{\dP}{\EM{\mathbb{P}}}
\newcommand{\dR}{\EM{\mathbb{R}}}
\newcommand{\dZ}{\EM{\mathbb{Z}}}

\newcommand{\al}{\alpha}
\newcommand{\ga}{\gamma}

\newcommand{\la}{\lambda}
\newcommand{\La}{\Lambda}
\newcommand{\si}{\sigma}
\newcommand{\de}{\delta}

\newcommand{\ABS}[1]{\EM{{\left| #1 \right|}}} 
\newcommand{\BRA}[1]{\EM{{\left\{#1\right\}}}} 
\newcommand{\NRM}[1]{\EM{{\left\| #1\right\|}}} 
\newcommand{\PAR}[1]{\EM{{\left(#1\right)}}} 
\newcommand{\SBRA}[1]{\EM{{\left[#1\right]}}} 

\begin{document}
\sloppy
\date{\today}

\title{H\"older regularity for operator scaling stable random fields}

\renewcommand{\thefootnote}{\fnsymbol{footnote}} \setcounter{footnote}{-1} 

\author{Hermine Bierm\'e}
\address{Hermine Bierm\'e, MAP5 Universit\'e Ren\'e Descartes, 45 rue des Saints P\`eres, 75270 Paris cedex 06, France}\email{hermine.bierme\@@{}math-info.univ-paris5.fr}
\urladdr{\url{http://www.math-info.univ-paris5.fr/~bierme/}}
\author{C\'eline Lacaux}
\address{C\'eline Lacaux, Institut \'Elie Cartan, UMR 7502, Nancy Universit\'e-CNRS-INRIA, BP 239, F-54506 Vandoeuvre-l\`es-Nancy, France}
\email{Celine.Lacaux\@@{}iecn.u-nancy.fr}
\urladdr{\url{http://www.iecn.u-nancy.fr/~lacaux/}}

\maketitle
\begin{abstract} We investigate the sample paths regularity of operator scaling $\al$-stable random fields. Such fields were introduced in~\cite{OSSRF} as anisotropic generalizations of self-similar fields and satisfy the scaling property $\{X(c^Ex) ; x\in\rd\}\dlim\{c^HX(x); x\in\rd\}$ where $E$ is a $d\times d$ real matrix
and $H>0$. In the case of harmonizable operator scaling random fields, the sample paths are locally H\"olderian and their H\"older regularity is characterized by the eigen decomposition of $\rd$ with respect to $E$. In particular, the directional H\"older regularity may vary and is given by the eigenvalues of $E$. In the case of moving average operator scaling random $\al$-stable random fields, with $\al\in(0,2)$ and $d\ge 2$, the sample paths are almost surely discontinous.   
\end{abstract}
\section{Introduction}
In this paper we consider operator scaling stable random fields as introduced in \cite{OSSRF}. 
More precisely, for $E$  a real $d\times d$ matrix with 
positive real parts of the eigenvalues, a scalar valued random field $\PAR{X(x)}_{x\in\rd}$ is called {\sl operator scaling} for $E$ and $H>0$ if for every $c>0$
\begin{equation}
\label{OSSRFintro}
\{X(c^Ex) ; x\in\rd\}\dlim\{c^HX(x); x\in\rd\},
\end{equation}
where $\dlim$ means equality of finite dimensional distributions and  as usual
 $c^E=\exp(E\log c)$ with 
$\exp(A)=\sum_{k=0}^\infty \frac{A^k}{k!}$ is the matrix exponential. Let us remark that up to consider the matrix $E/H$, we may assume, without loss of generality, that $H=1$.  

These
fields can be seen as anisotropic generalizations of self-similar random fields.
 Let us recall that a scalar valued random field $\PAR{X(x)}_{x\in\rd}$ is said to be {\it $H$-self-similar} with $H\in \dR$ if 
$$\{X(cx); x\in\rd\}\dlim\{c^HX(x); x\in\rd\}$$
for every $c>0$. Then a  $H$-self-similar field is also an operator scaling field for the matrix $E=I_d/H$, where $I_d$ is the identity matrix of size $d\times d$. 
Self-similar random fields are used in various  applications
such as internet traffic modelling \cite{wil:taqqu}, ground water modelling and
mathematical finance, just to mention a few. Various examples can be found for instance in the books \cite{levyvehel, Abry, sam:taqqu}.
A very important class of such fields are given by Gaussian random fields and especially by the fractional Brownian field~$B_H$,
where $H\in (0,1)$ is the so-called Hurst parameter. The random field $B_H$ is $H$-self-similar and  has stationary
increments, e.g. $\{B_H(x+h)-B_H(h);x\in\rd\}\dlim\{B_H(x);x\in\rd\}$  for any $h\in\rd$. It is an isotropic generalization of the famous fractional Brownian motion, implicitely introduced in \cite{Kolmo} and defined in \cite{MVN}.

Nevertheless, Gaussian random fields are not convenient for some heavy tails phenomena modelling. For this purpose, $\alpha$-stable random fields have been introduced. A scalar valued random field $\{X(x);x\in\rd\}$ is said to be symmetric $\alpha$-stable (S$\alpha$S), for $\alpha \in (0,2)$, if any linear combination $\underset{k=1}{\overset{n}{\sum}}a_kX(x_k)$ is S$\alpha$S. We address to the book \cite{sam:taqqu} for a well understanding
of such fields. The self-similar $\alpha$-stable fields with stationary increments have been extensively used to propose alternative to Gaussian modelling (see \cite{Mikosh,wil:taqqu} for instance) and are isotropic.

However,  certain applications (see, e.g., \cite{benson,Bonami} and references therein)
require that the random field is anisotropic and satisfies a
scaling relation. This scaling relation should have different
Hurst indices in different directions and these directions should
not necessarily be orthogonal. To reach this goal the authors of \cite{OSSRF} propose a new class of random fields with an anisotropic behavior driven by a $d\times d$ matrix $E$.  More precisely, they introduce  $\al$-stable random fields which have stationary increments and satisfy the operator scaling property~\eqref{OSSRFintro}. Two different classes of such fields are defined and analyzed, using a moving average
representation as well as an harmonizable
one. In the Gaussian case $\alpha=2$, according to \cite{OSSRF} there exist modifications of these fields which are
almost surely H\"older-continuous of certain indices. We give similar results here in the stable case $\alpha\in (0,2)$
for   harmonizable operator scaling stable random fields. Actually, we obtain their critical global and directional H\"older exponents, which are given by the eigenvalues of $E$. In general, such fields are anisotropic and their sample paths properties varies with the direction. In particular, in the case where $E$ is diagonalizable, for any eigenvector $\theta_j$ associated with the real eigenvalue $\la_j$, harmonizable operator scaling stable random field admits $H_j=1/\la_j$ as critical H\"older exponent in direction $\theta_j$. 
Let us point out that we establish an accurrate upper bound for the modulus of continuity.  Such upper bound has already been given in the case of real harmonizable fractional stable motions, which are self-similar, in~\cite{Kono} and in the case of some Gaussian random processes  in~\cite{Konog}. Then, in this paper, we generalize these results to any dimension $d$ and any harmonizable operator scaling stable fields. From this upper bound, we can deduce the H\"older sample paths regularity. Let us point out that we also  obtain this upper bound in the case of Gaussian operator scaling random field, which improves the sample paths properties establishes in~\cite{OSSRF}.

 Furthermore, whereas in the Gaussian case $\al=2$, moving average and harmonizable fields have the same kind of regularity properties, this is no more true in the case $\al\in(0,2)$. In particular, we show that for $d\ge 2$, a moving average operator scaling stable random field does not admit any continuous modification. Remark that if $d=1$, the sample paths regularity properties are already known since the processes studied  are self-similar moving average stable processes, see for example~\cite{sam:taqqu,Kono,Takashima}. \\

One of the main tool  for the study of sample paths of operator scaling random fields is the change of polar coordinates with respect to the matrix $E$
 introduced in~\cite{thebook} and used in \cite{OSSRF}. 
For $X$  a Gaussian  operator scaling random field   with stationary increments, using \eqref{OSSRFintro}, we can  write the variogramme of $X$ as
$$
v^2\PAR{x}=\dE\PAR{X^2\PAR{x}}=\tau_E\PAR{x}^{2H}\dE\PAR{X^2\PAR{\ell_E(x)}},
$$
where $\tau_E\PAR{x}$ is the radial part of $x$ with respect to $E$ and $\ell_E\PAR{x}$ is its polar part.  Therefore, in the Gaussian case, the sample paths regularity depends on the behavior of the polar coordinates $\PAR{\tau_E\PAR{x},\ell_E\PAR{x}}$ around $x=0$. Such property  also holds in the stable case $\al\in(0,2)$.  The H\"older regularity of the sample paths follows from estimates of $\tau_E\PAR{x}$ compared to $\NRM{x}$. These estimates are given in Section~\ref{PC} and their proof are postponed to the Appendix. 

Furthermore,  to study the sample paths in the stable case, the other main tool we use is a series representation of harmonizable operator scaling $\al$-stable random fields. Such representations in series of infinitely divisible laws have been studied in~\cite{LePage1,LePage2,Rosinskiu,Rosinskis}. As in~\cite{Kono}, our study  is based on a LePage series representation. Actually, the main idea is to choose a representation which is a conditionnally Gaussian series. \\

In Section~\ref{HR}, we recall the definition of harmonizable operator scaling random fields. Then, Sections~\ref{PC} and~\ref{LePagese} are devoted to the main tools we need for the study of the sample paths of these fields. More precisely, Section~\ref{PC} deals with the polar coordinates with respect to a matrix $E$ and Section~\ref{LePagese} gives the LePage series representation. In Section~\ref{HolderSE}, the sample paths properties and the Hausdorff dimension of the graph are given. Section~\ref{MASE} is concerned with moving average operator scaling random fields.
\section{Harmonizable representation}
\label{HR}
Let us recall the  definition of  harmonizable  operator scaling stable random fields, given
by \cite{OSSRF}.

Let $E$ be a real $d\times d$ matrix. Let $\la_1,\ldots,\la_d$ be the complex eigenvalues of $E$ and $a_j=\Re\PAR{\la_j}$ for each $j=1,\ldots,d$. We assume that 
\begin{equation}\label{VAP}
\displaystyle \min_{1\le j\le d} a_j>1.
\end{equation}
 
Let $\psi : \rd\rightarrow [0,\infty)$ be a continuous, $E^t$-homogeneous
function,  which means according to Definition~2.6 of \cite{OSSRF} that
$$\psi(c^{E^t}x)=c\psi(x) \mbox{ for all } c>0.$$
We assume moreover that  $\psi(x)\neq 0$ for $x\neq 0$. Such functions were studied in detail in~\cite{thebook}, Chapter~5 and various examples  are given in Theorem~2.11 and Corollary~2.12 of~\cite{OSSRF}.

Let $0<\alpha\le 2$ and $W_\al\PAR{d\xi}$ be a complex isotropic $\al$-stable random measure on $\dR^d$ with Lebesgue control measure (see~\cite{sam:taqqu} p.281). If $\al=2$, $W_\al\PAR{d\xi}$ is a complex isotropic Gaussian random measure. Let $q=\mbox{trace}\left(E\right)$. 

\begin{defi} Since \eqref{VAP} is fulfilled, the random field
\begin{equation}\label{harmonizable}
X_{\alpha}(x)=\Re\int_\rd\Bigl(e^{i\langle x,\xi
\rangle  }-1\Bigr)\psi(\xi)^{-1-q/\alpha}\,W_\alpha(d\xi)\ ,x\in\rd,
\end{equation}
is well defined and called harmonizable operator scaling stable random field.
\end{defi}
From~Theorem~4.1 and  Corollary~4.2 of~\cite{OSSRF}, $X_{\alpha}$ is stochastically continuous, has stationary increments and satisfies the following operator scaling property
\begin{equation}\label{OSP}
\forall c>0,\,\BRA{X_{\alpha}(c^Ex) ; x\in\rd}\dlim\BRA{cX_{\alpha}(x) ;x\in\rd}.
\end{equation}
 For notational sake of simplicity we denote the kernel function by
\begin{equation}\label{noyau}
f\PAR{x,\xi}=\PAR{\textup{e}^{i\langle x,\xi\rangle}-1}\psi\PAR{\xi}^{-1-q/\al}.
\end{equation}
Let us recall that $f\PAR{x,\cdot}\in L^{\alpha}(\rd)$ for any $x\in \rd$, which is a necessary and sufficient condition for~$X_{\alpha}$ to be defined.\\

Now, let us give some examples of such random fields.
\begin{exple}\label{expleK}{\rm
Let $I_d$ be the identity matrix of size $d\times d$, $E=I_d/H$ (with $0<H<1$) and $\psi\PAR{x}=\NRM{x}^H$ with $\|\cdot\|$ the Euclidean norm. Then the random field defined by \eqref{harmonizable} is a real harmonizable stable random field (see~\cite{sam:taqqu} for details on such fields). In this case, $X_{\alpha}$ is self-similar with exponent~$H$, i.e. 
$$\forall c>0,\,\{X_\al(cx); x\in\rd\}\dlim\{c^HX_\al(x); x\in\rd\}.$$
Let us quote that, if $\al=2$, $X_{\alpha}$ is 
a fractional Brownian field and its critical H\"older exponent is given by its Hurst index $H$  (see Theorem 8.3.2 of \cite{Adler} for instance). }
\end{exple}

\begin{exple}\label{diago}{\rm Assume that $E$ is diagonalizable. Then, all the eigenvalues of $E$ are real given by~$\left(a_j\right)_{1\le j\le d}$. Let $(\theta_j)_{1\le j\le d}$ be a basis of some corresponding eigenvectors and  consider the function $\psi$ defined by 
 $$\psi(x)=\left(\sum_{j=1}^d\left|\langle x,\theta_j\rangle\right|^{2/a_j}\right)^{1/2},\ x\in\rd.$$
 The function $\psi$ is clearly continuous and non negative on $\rd$. Moreover, since $\langle c^{E^t}x,\theta_j\rangle=\langle x,c^{E}\theta_j\rangle =c^{a_j}\langle x,\theta_j\rangle$, it is also $E^t$-homogeneous. Finally, since $(\theta_j)_{1\le j\le d}$ is a basis of $\rd$
 we have that $\psi(x)=0$ if and only if $x=0$. Then we can define $X_{\alpha}$ by \eqref{harmonizable}
and in this case the operator scaling property \eqref{OSP} implies that
$$\forall j=1,\ldots, d,\ \forall c>0, \ \BRA{X_\al\PAR{ct\theta_j} ; t\in\dR}\dlim \BRA{c^{1/a_j}X_\al\PAR{t\theta_j} ; t\in\dR}.
$$
The random field $X_\al$ is self-similar with Hurst index $H_j=1/a_j$ in the direction $\theta_j$. In particular, in the Gaussian case ($\al=2$), $\PAR{X_2\PAR{t\theta_j}}_{t\in\dR}$ is a fractional Brownian motion with Hurst index $H_j$. Then, in this case, its critical H\"older exponent is equal to $H_j$. 
}
\end{exple}

The main tool in the study of operator scaling random fields is the change of coordinates in a kind of polar coordinates
with respect to the matrix $E$. Then, before we study the sample paths of $X_{\alpha}$, we recall in the next section the definition of these coordinates and give some estimates of the radial part.

\section{Polar coordinates}
\label{PC}

According to Chapter~6 of~\cite{thebook},  
since $E$ is a real $d\times d$ matrix with positive real parts of the eigenvalues there exists
a norm $\|\cdot\|_E$ on $\rd$ such that for the unit sphere
$S_E=\{x\in\rd : \|x\|_E=1\}$ the map 
$$
\begin{array}{rccl}
 \Psi_E : &(0,\infty)\times
S_E&\longrightarrow & \rd\setminus\{0\}\\&(r,\theta)&\longmapsto&r^E\theta
\end{array}
$$
 is a homeomorphism.
Hence we can write any $x\in\rd\backslash\{0\}$
uniquely as $x=\tau_E(x)^E\ell_E(x)$ for some {\it radial part}
$\tau_E(x)>0$ and some direction $\ell_E(x)\in S_E$ such that
$x\mapsto\tau_E(x)$ and $x\mapsto \ell_E(x)$ are continuous. Observe that
$S_E=\{x\in\rd : \tau_E(x)=1\}$ is compact and set 
\begin{equation}\label{ME}
m_E=\underset{S_E}{\min}\|x\| \ \mbox{ and } \ M_E=\underset{S_E}{\max}\|x\|.
\end{equation}
We know that
$\tau_E(x)\to\infty$ as $x\to\infty$ and $\tau_E(x)\to 0$ as $x\to 0$.
Hence we can extend $\tau_E(\cdot)$ continuously by setting
$\tau_E(0)=0$. \\
Let us recall the formula of integration in {\it polar coordinates} established in \cite{OSSRF}. 
\begin{prop}\label{CDV}
There exists a unique finite Radon measure $\sigma_E$ on $S_E$ such
that for all $f\in L^1(\rd,dx)$ we have
\[\int_\rd f(x)\,dx=\int_0^\infty\int_{S_E}f(r^E\theta)\,\sigma_E(d\theta)\,r^{q-1}\,dr
.\]\\
\end{prop}

In the Gaussian case ($\al=2$), the variogramme of $X_2$ can be rewritten as follows 
$$
v^2\PAR{x}=\tau_E(x)^2\mathbb{E}\left(X_2\left(\ell_E(x)\right)^2\right),
$$
using the operator scaling property~\eqref{OSP}. 
Then, the H\"older regularity of $X_2$ follows from estimates of~$\tau_E\PAR{x}$ compared to $\NRM{x}$ around $x=0$, e.g. the H\"older regularity of $\tau_E$ around $0$, see~\cite{OSSRF}. Then, in order to get an upper bound for the modulus of continuity (for any $\al$), we need  precise estimates of~$\tau_E\PAR{x}$. 

As done in \cite{Xiao} for the study of operator-self-similar Gaussian
 random fields we use the Jordan decomposition of the matrix $E$ to get estimates of $\tau_E$. From the Jordan decomposition's theorem (see \cite{Hirsch} p. 129 for instance), 
there exists a real invertible~$d\times d$ matrix $P$ such that $D=P^{-1}EP$ is
of the real canonical form, which means that $D$ is 
 composed of diagonal blocks which are either Jordan cell matrix of the form
 $$
\left(\begin{array}{ccccc}
\lambda &1&\ldots& 0\\
0&\lambda&  \ddots& \vdots\\
\vdots&\ddots& \ddots&1\\
0&\ldots&0&\lambda
\end{array}\right)
$$ with $\lambda$ a real eigenvalue of $E$ or blocks of the form
\begin{equation}
\label{Jordan block}
\left(\begin{array}{cccccc}
\Lambda &I_2&0&\ldots& 0\\
0&\Lambda& I_2& \ddots& \vdots\\
\vdots & \ddots& \ddots& \ddots& 0\\
\vdots&\ddots&\ddots & \ddots&I_2\\
0&\ldots&\ldots&0&\Lambda
\end{array}\right) \mbox{ with }
\Lambda=\left(\begin{array}{cc} a & -b\\
b & a\end{array}\right) \mbox{ and } I_2=\left(\begin{array}{cc} 1 & 0\\
0 & 1\end{array}\right),
\end{equation}
where the complex numbers $a \pm i b$ ($b\neq 0$) are complex conjugated eigenvalues of $E$. \\

Let us denote by $\|\cdot\|$ the subordinated norm of the Euclidean norm on the matrix space. 
Precise estimates of $\tau_E$ follow from the next lemma.
\begin{lemma}\label{Jordan power}
Let $J$ be either a Jordan cell matrix of  size $l$ or a block of the form~\eqref{Jordan block} of size $2l$ associated with the eigenvalue $\lambda=a+ib$. 
Then, for any $t\in (0,e^{-1}]\cup[e,+\infty)$
$$t^{a}\le \|t^{J}\|\le \sqrt{2 l}\,e\, t^{a}\left|\log t\right|^{l-1}.$$
\end{lemma}
\begin{proof} see the Appendix.
\end{proof}
 Let us be more precise on the Jordan decomposition of $E$. Let us recall that the eigenvalues of $E$ are denoted by $\la_j$, $j=1\ldots d$ and that $a_j=\Re\PAR{\la_j}>1$ for $j=1,\ldots,d$. There exist $J_1,\ldots,J_m$, where each $J_j$ is either a Jordan cell matrix or a block of the form~\eqref{Jordan block}, and $P$ a real $d\times d$ invertible matrix such that 
 $$
 E=P\left(\begin{array}{cccc}J_1 & 0 &\ldots & 0\\
0&J_2& 0&\vdots\\
\vdots&\ddots& \ddots& 0\\
0&\ldots& 0& J_m\end{array}\right)P^{-1}.
 $$
We can assume  that each $J_j$ is associated with the eigenvalue $\la_j$ of $E$ and that 
 $$
1<a_1\le \cdots\le a_m.
$$
We also set $H_j=a_j^{-1}$ and have
\begin{equation}\label{Hj}
0<H_m\le\cdots\le H_1<1. 
\end{equation}
 If $\lambda_j\in\dR$, $J_j$ is a Jordan cell matrix of size $\tilde{l_j}=l_j\in \mathbb{N}\backslash\{0\}$. If $\lambda_j\in\dC\backslash\dR$, $J_j$ is a block of the form~\eqref{Jordan block} of size~$\tilde{l_j}=2l_j\in 2\mathbb{N}\backslash\{0\}$. 
Then for any $t>0$,
$$t^E=P\left(\begin{array}{cccc}t^{J_1} & 0 &\ldots & 0\\
0&t^{J_2}& 0&\vdots\\
\vdots&\ddots& \ddots& 0\\
0&\ldots& 0& t^{J_m}\end{array}\right)P^{-1}
$$ 
We denote by $\PAR{e_1,\ldots,e_d}$ the canonical basis of $\rd$ and set $f_j=Pe_j$ for every $j=1,\ldots,d$. Hence,~$\PAR{f_1,\ldots,f_d}$ is a basis of $\dR^d$. For all $j=1,\ldots,m$,  let  
\begin{equation}
\label{Wj}
W_j=\mbox{Vect}\left(f_k\,\,;\,\,\underset{i=1}{\overset{j-1}{\sum}}\tilde{l_{i}}+1\le k\le \underset{i=1}{\overset{j}{\sum}}\tilde{l_j}\right).
\end{equation}  Then,
each $W_j$ is a $E$-invariant set and $E=\bigoplus_{j=1}^m W_j$.\\

 The following result gives bounds on the growth
rate of $\tau_E(x)$ in terms of the real parts of the eigenvalues of
$E$.

\begin{prop}\label{taulocal} For any $r\in(0,1)$ there exist  $c_1, c_2>0$ such that for every $1\le j_0\le j\le m$, 
$$ c_1 \|x\|^{H_{j_0}}\left|\log \|x\|\right|^{-(p_{j_0,j}-1)H_{j_0}}\le \tau_E(x)\le c_2\|x\|^{H_j}\left|\log \|x\|\right|^{(p_{j_0,j}-1)H_j}
$$
holds for any  $x\in \oplus_{k=j_0}^{j}W_k$ with $\|x\|\le r$ and $p_{j_0,j}=\underset{j_0\le k\le j}{\max}l_k$.
\end{prop}

\begin{proof} See the Appendix.
\end{proof}

Then, we easily deduce the following corollary. 
\begin{cor}\label{tauglobal} For any $r\in(0,1)$ there exist  $c_1, c_2>0$ such that for any $x\in W_j$ with $\|x\|\le r$ 
\begin{equation}
\label{control WjE} 
c_1 \|x\|^{H_j}\left|\log \|x\|\right|^{-(l_j-1)H_j}\le \tau_E(x)\le c_2\|x\|^{H_j}\left|\log \|x\|\right|^{(l_j-1)H_j}
\end{equation} 
and for any $x\in \mathbb{R}^d$ with $\|x\|\le r$
\begin{equation}\label{control uE} c_1 \|x\|^{H_1}\left|\log \|x\|\right|^{-(l-1)H_1}\le \tau_E(x)\le c_2\|x\|^{H_m}\left|\log \|x\|\right|^{(l-1)H_m}
\end{equation} 
  where $l=\underset{1\le j\le m}{\max}l_j$.
\end{cor}

Therefore we have precise estimates for the H\"older regularity of the radial part. Let us remark that we improve 
the first statement of Lemma 2.1 of \cite{OSSRF} and that the second one can also be improved in a similar way.
 From these  estimates  we deduce  the H\"older regularity of $X_\al$ in Section~\ref{HolderSE}.  As already mentionned, the study of the sample paths is based on a  series representation. Then, before we state regularity properties, it remains to give the LePage series representation of harmonizable operator scaling random fields $X_\al$ defined by~\eqref{harmonizable}.

\section{Representation as a LePage series}
\label{LePagese}
An overview on series representations of infinitely divisible random variable without Gaussian part can be found for example  in~\cite{Rosinskis,Rosinskipp} and references therein.  In particular, LePage series representation (\cite{LePage1,LePage2}) have been used in~\cite{Konot,Kono} to study the sample paths regularity of some self-similar $\al$-stable random motions with $\al\in(0,2)$.  Here, this representation is also the main representation we use in the case $\al\in(0,2)$. Actually, in the Gaussian case $\al=2$, such representation does not hold.\\

Let us now introduce some notations that will be used throughout the paper. Let $\mu$ be an arbitrary probability measure equivalent to the Lebesgue measure on $\rd$ and let $m$ be its Radon-Nikodym derivative that is $\mu\PAR{d\xi}=m\PAR{\xi}d\xi$.

\noindent {\bf Notation } Let ${\PAR{T_n}}_{n\ge 1}$, 
${\PAR{g_n}}_{n \ge 1}$ and ${\PAR{\xi_n}}_{n \ge 1}$ be independent sequences of random variables. 
\begin{enumerate}
\item[$\bullet$] $T_n$ is the $n$th arrival time of a Poisson process with  intensity 1.
\item[$\bullet$] ${\PAR{g_n}}_{n \ge 1}$ is a sequence of i.i.d. isotropic complex random variables so that 
$g_n\stackrel{(d)}{=}\mbox{e}^{i\theta}g_n$ for any $\theta\in\dR$. We also assume that $0<\dE\PAR{\ABS{g_n}^\alpha}<+\infty$.
\item[$\bullet$] ${\PAR{\xi_n}}_{n \ge 1}$ is a sequence of i.i.d. random
  variables with common law $  \mu\PAR{d\xi}=m\PAR{\xi}d\xi. $\\
\end{enumerate}

According to Chapter 3 and Chapter 4 of~\cite{sam:taqqu}, stochastic integrals with respect to an $\al$-stable  random measure $\La$ can be represented as a LePage series as soon as the control measure of $\La$ is a finite measure. The next proposition generalizes this to $W_\al$ whose control measure is the Lebesgue measure. It is a consequence of Lemma~4.1 of~\cite{Konot}, which is a correction of Lemma 1.4 of \cite{Marcus}. This proposition can also be deduced from~\cite{Rosinskis,Rosinskiu}, concerned with  series representations of stochastic integrals with respect to infinitely divisible random measures.   
\begin{prop}\label{Lepage1}
Assume that $\al\in(0,2)$. Then, for every complex valued function $h\in L^\al\PAR{\dR^d}$, the series 
\begin{equation}
Y^h=\sum_{n=1}^{+\infty}T_n^{-1/\al}m\PAR{\xi_n}^{-1/\al}h\PAR{\xi_n}g_n   
\end{equation}
converges almost surely. Furthermore, 
$$
C_\al Y^h\stackrel{\PAR{d}}{=}\int_{\dR^d}\!\!h\PAR{\xi}W_\al\PAR{d\xi}
$$
with 
\begin{equation}
\label{defcste}
C_\al=\dE\PAR{\ABS{\Re\PAR{g_1}}^\al}^{-1/\al}\PAR{\frac{1}{2\pi}\int_{0}^\pi\ABS{\cos\PAR{x}}^\al dx}^{1/\al}\PAR{\int_{0}^{+\infty}\frac{\sin\PAR{x}}{x^\al}dx}^{-1/\al}. 
\end{equation}
\end{prop}
\begin{rem} According to Proposition~\ref{Lepage1}, taking $\al\in(0,2)$,  the random measure 
$$
\La_\al\PAR{d\xi}=C_\al\sum_{n=1}^{+\infty}T_n^{-1/\al}m\PAR{\xi_n}^{-1/\al}g_n\de_{\xi_n}\PAR{d\xi}
$$
is a complex isotropic $\al$-stable random measure. 
\end{rem}
\begin{proof} 
Let $V_n=m\PAR{\xi_n}^{-1/\al}h\PAR{\xi_n}g_n$. Then, $V_n$, $n\ge 1$, are i.i.d. isotropic complex random variables. By Lemma~4.1 in~\cite{Konot}, $Y^h$ converges almost surely and 
$$
\forall z\in\dC,\, \dE\PAR{\exp\PAR{i\Re\PAR{\bar{z}Y^h}}}=\exp\PAR{-\sigma^\al\ABS{z}^\al}
$$
with 
$$
\si^\al=\dE\PAR{\ABS{\Re\PAR{V_1}}^\al}\int_{0}^{+\infty}\frac{\sin\PAR{x}}{x^\al}dx.
$$
Since $g_1$ is invariant by rotation and independent with $\xi_1$,
$$
\dE\PAR{\ABS{\Re\PAR{V_1}}^\al}=\dE\PAR{m\PAR{\xi_1}^{-1}\ABS{h\PAR{\xi_1}}^\al}\dE\PAR{\ABS{\Re\PAR{g_1}}^\al}=\dE\PAR{\ABS{\Re\PAR{g_1}}^\al}\int_{\dR^d}\ABS{h\PAR{\xi}}^\al d\xi.
$$
Moreover, by definition of an isotropic $\al$-stable random stable measure (see \cite{sam:taqqu}),
$$
\forall z\in\dC, \dE\PAR{\exp\PAR{i\Re\PAR{\bar{z}\int_{\dR^d}h\PAR{\xi}M\PAR{d\xi}}}}=\exp\PAR{-c_\al^\al\PAR{h}\ABS{z}^\al}
$$
with $c_\al^\al\PAR{h}=\PAR{\displaystyle\frac{1}{2\pi}\int_{0}^\pi\ABS{\cos\PAR{x}}^\al dx }\int_{\dR^d}\ABS{h\PAR{\xi}}^\al d\xi.$  Therefore, we have 
$$
C_\al Y^h\stackrel{\PAR{d}}{=}\int_{\dR^d}\!\!h\PAR{\xi}W_\al\PAR{d\xi}
$$
where $C_\al$ is defined by~\eqref{defcste}, which concludes the proof.
\end{proof}

From the previous proposition, we deduce the following statement which is the main series representation we use in our investigation.
\begin{prop}
\label{LePage} Let $\al\in (0,2)$. 
For every $x\in\dR^d$, the series 
\begin{equation}
\label{GNSN}
Y_\al(x)=C_\al\Re\PAR{\sum_{n=1}^{+\infty}T_n^{-1/\al}m\PAR{\xi_n}^{-1/\al}f\PAR{x,\xi_n}g_n},
\end{equation}
where $C_\al$ is defined by~\eqref{defcste}, 
converges almost surely. Furthermore, 
$$
\BRA{Y_\al(x): x\in \dR^d}\dlim \BRA{X_\al(x): x\in \dR^d}.
$$
\end{prop}
\begin{proof} From Proposition \ref{Lepage1}, for any $x\in \rd$, the convergence of the series follows from
the fact that $f(x,\cdot)\in  L^\al\PAR{\dR^d}$. The equality in distribution between $X_\al$ and $Y_\al$ is obtained by linearity of both fields.
\end{proof}

Using LePage representation~\eqref{GNSN} of $X_\al$ and  the estimates given in Section~\ref{PC}, we  give an upper bound for the modulus of continuity of $X_\al$ and obtain the critical H\"older regularity of its sample paths in the next section.

\section{H\"older regularity and Hausdorff dimension}\label{HolderSE}
Throughout this section we fix $K$ a compact set of $\rd$ and investigate the H\"older regularity of the sample paths of $X_{\alpha}$ on $K$, with $X_{\alpha}$ the harmonizable operator scaling stable random field  defined by \eqref{harmonizable}.

Let us recall that for the Gaussian case $\alpha=2$, according to Theorem 5.4 of \cite{OSSRF}, the H\"older
regularity of $X_2$ depends on the subspaces $\left(W_j\right)_{1\le j\le m}$ defined by 
\eqref{Wj} and associated to the eigenvalues of~$E$.
 More precisely, Theorem 5.4 of \cite{OSSRF} implies that,
when restricted to the subspace~$W_j$, the Gaussian random field $\left\{X_2(x); x\in W_j\right\}$ admits $H_j$ as critical H\"older exponent. This follows from the fact that the regularity of $X_2$ is determined by the regularity of $\tau_E$ around $0$, which is given by $H_j$ according to~\eqref{control WjE}. Here, we give an upper bound for the modulus of continuity of
$X_{\alpha}$  in the general case~$\alpha\in (0,2]$. Then we prove that the critical H\"older exponents
are the same than in the Gaussian case~$\alpha=2$.
Let us state our main result.
\begin{theorem}\label{ModCont} Let $\alpha\in (0,2)$.
There exists  a modification $X^*_{\alpha}$ of $X_\al$ on $K$ such that 
\begin{equation}\label{module}
\lim_{\delta\downarrow 0}\underset{\underset{\|x-y\|\le \delta}{x,y\in K}}{\sup}\frac{\ABS{X^*_\al(x)-X^*_\al(y)}}{\tau_E(x-y)\ABS{\log \tau_E(x-y)}^{1/\alpha+1/2+\varepsilon}}=0  \mbox{ a.s.}
\end{equation}
for any $\varepsilon>0$. 
\end{theorem}
This result was proved in the case of harmonizable self-similar processes in~\cite{Kono}, e.g. in the case of Example~\ref{expleK} with $d=1$. The main idea is to use the LePage series representation~\eqref{GNSN} where~$g_n$, $n\ge 1$, are Gaussian  complex isotropic random variables. Furthermore, it remains to choose the density distribution $m$ of the $\xi_n$. In~\cite{Kono}, the authors choose
$$
m\PAR{\xi}=\frac{c_\eta}{\ABS{\xi}\ABS{\log\ABS{\xi}}^{1+\eta}}, \ \xi\in \mathbb{R}\backslash\{0\}
$$  
where $c_\eta>0$. A straightforward generalization in higher dimension $d$ leads to choose 
$$
m\PAR{\xi}=\frac{c_\eta}{\NRM{\xi}^d\ABS{\log\NRM{\xi}}^{1+\eta}}, \ \xi\in\rd\backslash\{0\}. 
$$ 
Remark that in this case (e.g. Example~\ref{expleK}) the matrix $E=I_d/H=E^t$ and that we can choose $\NRM{\,\cdot\,}_{E^t}=\NRM{\,\cdot\,}$. Then, using classical polar coordinates, we obtain that for $x\ne 0$, 
$$\tau_{E^t}\PAR{x}=\NRM{x}^{H} \ \textrm{and} \ \ell_{E^t}\PAR{x}=\frac{x}{\NRM{x}}$$
and therefore 
$$
m\PAR{\xi}=\frac{c_\eta}{\tau_{E^t}\PAR{\xi}^q\ABS{\log\tau_{E^t}\PAR{\xi}}^{1+\eta}}
$$
since $q={\rm{trace}}\PAR{E}=d/H$. 
Then we choose this density in the general case. The advantage is that~$m$ only depends on the radial part $\tau_{E^t}$. 
\begin{proof}[Proof of Theorem \ref{ModCont}] We can assume without loss of generality that $K=[0,1]^d$. 
According to Proposition~\ref{LePage} for every $x\in\rd$
$$
Y_\al\PAR{x}=C_\al\Re\PAR{\sum_{n=1}^{+\infty} T_n^{-1/\al}m\PAR{\xi_n}^{-1/\al}f\PAR{x,\xi_n}g_n}
$$
converges almost surely and  $Y_\al \dlim X_\al$. As already mentionned, we assume that $g_n$, $n\ge 1$ are  Gaussian complex isotropic random variables and that the density distribution of $\xi_n$ is defined by  
$$
m\PAR{\xi}=\frac{c_\eta}{\tau_{E^t}\PAR{\xi}^q\ABS{\log\tau_{E^t}\PAR{\xi}}^{1+\eta}},\ \xi\in\rd\backslash\{0\},
$$
where $\eta>0$ and $c_\eta$ is chosen such that $\int_{\rd}m\PAR{\xi}d\xi=1$. In particular, conditionally to $\PAR{T_n,\xi_n}_n$,~$Y_{\alpha}(x)$ is a real-valued Gaussian random variable.\\

As in the proof of Kolmogorov-Centsov Theorem (see~\cite{karatzas}), we first prove that almost surely 
for~$\tau_E\PAR{x_k-x_{k'}}$ small enough, 
$$
\left|Y_\al(x_k)-Y_\al(x_{k'})\right|\le C \tau_E(x_k-x_{k'})\ABS{\log \tau_E(x_k-x_{k'})}^{1/\alpha+1/2+\varepsilon}.
$$ 
where $\PAR{x_k}_k$ is some countable dense sequence  of $K$. Then, $X_\al$ satisfies the same property. Finally, we give a modification $X_\al^*$ of $X_\al$ for which~\eqref{module} holds. In the first step, we construct the sequence~$\PAR{x_k}_k$ we use.  \\

{\bf{Step 1.}} 
For $k\in\mathbb{N}\backslash\{0\}$ let us choose $\nu_k\in \mathbb{N}\backslash\{0\}$ the smallest integer such that
$$
c_2 d^{H_m/2}2^{-\nu_kH_m}\ABS{\nu_k\log 2}^{(l-1)H_m}\le 2 ^{-k},
$$
where $c_2$ and $l$  are given by  Corollary \ref{tauglobal}
 for $r=1/2$. Remark that by definition, $\nu_k\le \nu_{k+1}$. Up to change $c_2$ in Proposition~\ref{tauglobal},  we can assume that 
$$
c_2 d^{H_m/2}2^{-H_m}\ABS{\log 2}^{(l-1)H_m}>1,
$$
which implies that $ \nu_k>1$ for every $k$. Then, since $1\le \nu_k-1\le \nu_k$,
$$
c_2 d^{H_m/2}2^{-\PAR{\nu_k-1}H_m}\ABS{\PAR{\nu_k-1}\log 2}^{(l-1)H_m}> 2 ^{-k}
$$
and  we have 
\begin{equation}
\label{nuk}
2^{-k}{\left(2\sqrt{d}\right)^{-H_m}c_2^{-1}}<\left(2^{-\nu_k}\left(\nu_k\log 2\right)^{l-1}\right)^{H_m}\le 
2^{-k}{\left(\sqrt{d}\right)^{-H_m}c_2^{-1}}.
\end{equation}
Let us remark that $\PAR{\nu_k}_{k\ge 1}$ is an increasing sequence and then that $\nu_k\ge k$ for every $k$. Furthermore, taking the logarithm of~\eqref{nuk}, one easily proves that 
$$
\lim_{k\to +\infty}\frac{k}{\nu_k}=H_m, 
$$
which implies $\nu_k=O\PAR{k}$. 

For every $k\in  \dN\backslash\{0\}$ and $j=(j_1,\ldots,j_d)\in \mathbb{Z}^d$ 
we set 
$$
x_{k,j}=\frac{j}{2^{\nu_k}}, \ {\mathcal D}_k=\displaystyle\BRA{x_{k,j} :  j\in\dZ^d\cap\SBRA{0,2^{\nu_k}}^d} \ \textrm{and} \ N_k=\rm{card}\,\mathcal{D}_k= \PAR{2^{\nu_k}+1}^d.
$$ 
Then $\mathcal{D}_k$ is a $2^{-k}$ net of $K$ for $\tau_E$ in the sense that for any $x\in K$ one can find $x_{k,j}\in\mathcal{D}_k$ such that $\tau_E(x-x_{k,j})\le 2^{-k}$. Actually, by Proposition~\ref{tauglobal}, it is sufficient to  choose $j$ such that $j_i\le 2^{\nu_k}x_i<j_i+1$ for $1\le i\le d$. 

Let us remark that the sequence $\PAR{{\mathcal D}_k}_k$ is increasing and set ${\mathcal D}=\displaystyle\bigcup_{k=1}^{+\infty}{\mathcal D}_k$.

{\bf{Step 2.}}
Almost surely, for any $x,y\in\cD$
$$
{Y}_\al\PAR{x}-{Y}_\al\PAR{y}=C_\al\Re\PAR{\sum_{n=1}^{+\infty}T_n^{-1/\al}m\PAR{\xi_n}^{-1/\al}\PAR{f\PAR{x,\xi_n}-f\PAR{y,\xi_n}}g_n},
$$
where $C_\al$ is defined by \eqref{defcste}.
Since the sequences $\PAR{T_n}_n$, $\PAR{\xi_n}_n$ and $\PAR{g_n}_n$ are independent and $\PAR{g_n}_n$
is a sequence of i.i.d. Gaussian complex isotropic random variables
$$
R(x,y)=\sum_{n=1}^{+\infty}T_n^{-1/\al}m\PAR{\xi_n}^{-1/\al}\PAR{f\PAR{x,\xi_n}-f\PAR{y,\xi_n}}g_n
$$
is a Gaussian isotropic complex random variable conditionally to $\PAR{T_n,\xi_n}_n$. 
Remark that ${Y}_\al\PAR{x}-{Y}_\al\PAR{y}=C_\al\Re\PAR{R(x,y)}$ almost surely.
Therefore, conditionally
to $\PAR{T_n,\xi_n}_n$, ${Y}_\al\PAR{x}-{Y}_\al\PAR{y}$ is a real centered Gaussian  random variable
with variance 
$$\begin{array}{rcl}
\displaystyle v^2\PAR{\PAR{x,y}\left|\right.\PAR{T_n,\xi_n}_n}&=&\displaystyle \frac{C_\al^2}{2}\,\dE\PAR{\ABS{R\PAR{x,y}}^2\left|\right.\PAR{T_n,\xi_n}_n}\\
&=&\displaystyle\frac{C_\al^2}{2}\,\dE\PAR{\ABS{g_1}^2}\sum_{n=1}^{+\infty}T_n^{-2/\al}m\PAR{\xi_n}^{-2/\al}\ABS{f\PAR{x-y,\xi_n}}^2,
\end{array}
$$
since $\ABS{f\PAR{x,\xi_n}-f\PAR{y,\xi_n}}=\ABS{f\PAR{x-y,\xi_n}}$.

 We  consider  the
set 
$$
E_{i,j}^k=\left\{\omega : \left|Y_\al\PAR{x_{k,i}}-Y_\al\PAR{x_{k,j}}\right|> v\PAR{\PAR{x_{k,i},x_{k,j}}\left|\right.\PAR{T_n,\xi_n}_n}\,\varphi\PAR{\tau_E\PAR{x_{k,i}-x_{k,j}}}
\right\},
$$
where, as in \cite{Konog}, 
$$
\varphi\PAR{t}=\sqrt{2C_\varphi d\log\frac{1}{t}}, \ t>0
$$
and $C_\varphi>0$ will be chosen later. Then, for every $(k,i,j)$, 
$$
\dP\PAR{E_{i,j}^k}=\dE\PAR{\dE\PAR{{\bf{1}}_{E_{i,j}^k}\left|\right.\PAR{T_n,\xi_n}_n}}.
$$
We give in the following an upper bound of this probability for well chosen $\PAR{k,i,j}$. Note that if $Z$ is a real centered Gaussian random variable with variance $1$, we have 
$$
\dE\PAR{{\bf{1}}_{E_{i,j}^k}\left|\right.\PAR{T_n,\xi_n}_n}=\dP\PAR{\ABS{Z}>\varphi\PAR{\tau_E\PAR{x_{k,i}-x_{k,j}}}} \ \mbox{almost surely}. 
$$
Let us choose $\delta \in \left(0,1\right)$ and 
set for $k\in \mathbb{N}\backslash\{0\}$, $\delta_k=2^{-(1-\delta)k}$ and
$$I_k=\BRA{\PAR{i,j}\in\dZ^d\cap\SBRA{0,2^{\nu_k}}^d :  \tau_E\PAR{x_{k,i}-x_{k,j}}\le \delta_k}.$$
For every $\PAR{i,j}\in I_k$, since $\varphi$ is a decreasing function
$$
\dP\PAR{\ABS{Z}>\varphi\PAR{\tau_E\PAR{x_{k,i}-x_{k,j}}}} \le \dP\PAR{\ABS{Z}>\varphi\PAR{\de_k}}. 
$$
 Then, we recall that 
$$
\forall u\ge 0, \dP\PAR{Z>u}\le\frac{\textup{e}^{-u^2/2}}{\sqrt{2\pi}u}.
$$
Therefore, for every $k\in\dN\backslash\{0 \}$ and $\PAR{i,j}\in I_k$,  
$$
\mathbb{P}\left(E_{i,j}^k \right)\le\sqrt{\frac{2}{\pi}}\frac{\textup{e}^{-\varphi^2\PAR{\delta_k}/2}}{\varphi\PAR{\delta_k}}
= \frac{2^{-\PAR{1-\de}kC_\varphi d}}{\sqrt{{2C_{\varphi}d\PAR{1-\de}k\log 2}}}.$$
Hence, 
$$
\sum_{k=1}^{\infty}\sum_{\PAR{i,j}\in I_k}\mathbb{P}\left(E_{i,j}^k \right)\le \frac{1}{\sqrt{2C_{\varphi}d\PAR{1-\de}\log 2}}\sum_{k=1}^{+\infty} 2^{-\PAR{1-\de}kC_\varphi d} \,\rm{card}\, I_k. 
$$
Since $\nu_k=O\PAR{k}$, the lower bounds of \eqref{nuk} and Corollary \ref{tauglobal} leads to
 $$\textrm{card}\left\{j\in \dZ^d\cap\SBRA{0,2^{\nu_k}}^d : \left(i,j\right)\in I_k\right\}=O\left(\delta_k^{d/H_1}2^{dk/H_m}k^{2d(l-1)}\right),$$
for any $ i\in \dZ^d\cap\SBRA{0,2^{\nu_k}}^d$. 
Then one can find a finite constant $C>0$ such that
$$\sum_{k=1}^{\infty}\sum_{\PAR{i,j}\in I_k}\mathbb{P}\left(E_{i,j}^k \right)\le 
\frac{C}{\sqrt{C_{\varphi}\PAR{1-\de}}}\sum_{k=1}^{\infty}{k^{3d(l-1)}}2^{-kd\left(-\frac{2}{H_m}+\frac{1-\delta}{H_1}+(1-\delta)C_{\varphi} \right)},
$$
which is finite 
for $C_{\varphi}> \frac{2}{H_m}-\frac{1}{H_1}$ and $\delta$ small enough.
 By the Borel-Cantelli Lemma, almost surely there exists an integer $ k^*\PAR{\omega}$  such that for every $k\ge k^*\PAR{\omega}$,
$$
\ABS{Y_\al\PAR{x}-Y_\al\PAR{y}}\le v\PAR{\PAR{x,y}\left|\right.\PAR{T_n,\xi_n}_n}\,\varphi\PAR{ \tau_E\PAR{x-y}}
$$
as soon as $x, y \in{\mathcal D}_k$ with $\tau_E\PAR{x-y}\le \de_k$.  \\
\\

{\bf Step 3.}
As in \cite{Kono} let us give an upper bound of 
$$v^2\PAR{\PAR{x,y}\left|\right.\PAR{T_n,\xi_n}_n}=\displaystyle\frac{C_\al^2}{2}\dE\PAR{\ABS{g_1}^2}\sum_{n=1}^{+\infty}T_n^{-2/\al}m\PAR{\xi_n}^{-2/\al}\ABS{f\PAR{x-y,\xi_n}}^2$$ with respect to $\tau_E\PAR{x-y}$.
By definition of $f$
$$v^2\PAR{\PAR{x,y}\left|\right.\PAR{T_n,\xi_n}_n}\le \displaystyle\frac{C_\al^2}{2}\dE\PAR{\ABS{g_1}^2}\sigma^2\left(\tau_E\PAR{x-y}\right),$$
where for all $h>0$
$$\sigma^2(h)=\sum_{n=1}^{+\infty}T_n^{-2/\al}m\PAR{\xi_n}^{-2/\al}\min\PAR{M_E\left\|h^{E^t}\xi_n\right\|,2}^2\psi\PAR{\xi}^{-2-2q/\al},$$
 with $M_E$ given by \eqref{ME}. For the sake of clearness we postpone the proof of the control
of $\sigma^2(h)$ in Appendix and state it in the following lemma.
\begin{lemma}\label{esperance} For any $\gamma\in(0,1)$ there exists a finite constant $c>0$ such that 
 $$
  \mathbb{E}\PAR{\sigma^2(h)\left|\right.\PAR{T_n}_n}\le c\sum_{n=1}^{+\infty}T_n^{-2/\al} h^{2}\ABS{\log h}^{\PAR{1+\eta}\PAR{2/\al-1}} \mbox{ almost surely }
 $$
as soon as $h\le 1-\gamma$.
\end{lemma}
Following  \cite{Kono} let us denote
$$b(h)=h|\log h|^{(1+\eta)/\alpha}.$$
Then by Lemma \ref{esperance}, 
$$\mathbb{E}\PAR{\left.\sum_{k=1}^{+\infty}\frac{\sigma^2(2^{-k})}{b^2(2^{-k})}\right|\PAR{T_n}_n}<+\infty.$$ Therefore by independence of $\PAR{T_n}_n$ and $\PAR{\xi_n}_n$, almost surely 
$$
 \lim_{k\to +\infty}\frac{\sigma^2(2^{-k})}{b^2(2^{-k})}=0.
$$
Up to change the Euclidean norm $\NRM{\,\cdot\,}$ by the equivalent norm $\NRM{\,\cdot\,}_{E^t}$ defined in Lemma~6.1.5 
of~\cite{thebook} the map $h\mapsto \NRM{h^{E^t}\xi}$ is increasing and so is $\si^2$. Also, one can conclude, as  in~\cite{Kono}, that almost surely
$$\lim_{h\rightarrow 0}\frac{\sigma^2(h)}{b^2(h)}=0.$$
Therefore, up to change  $k^*$ one can assume that for every $k\ge k^*\PAR{\omega}$, for every $x,y\in\cD_k$,
\begin{equation}\label{majoY}
\ABS{Y_\al\PAR{x}-Y_\al\PAR{y}}\le \sqrt{2dC_{\varphi}}\tau_E\PAR{x-y}\ABS{\log \tau_E\PAR{x-y}}^{(1+\eta)/\alpha+1/2}.
\end{equation}
as soon as $\tau_E\PAR{x-y}\le \de_k$. Let
$$
\Omega^*=\bigcup_{n=1}^{+\infty}\bigcap_{k\ge n}\bigcap_{\underset{\tau_E\PAR{x-y}\le \de_{k}}{x,y\in\cD_{k}} }\BRA{\ABS{X_\al\PAR{x}-X_\al\PAR{y}}\le \sqrt{2dC_{\varphi}}\tau_E\PAR{x-y}\ABS{\log \tau_E\PAR{x-y}}^{(1+\eta)/\alpha+1/2}}
$$ 
Since $X_\al$ and $Y_\al$ have the same  finite dimensional margins $\mathbb{P}\PAR{\Omega^*}=1$.   \\
\\

{\bf{Step 4.}} Let $\omega\in \Omega^*$. By Step 3 there exists $k^*(\omega)\ge 1$ such that $X_\al$
satisfies \eqref{majoY} for $k\ge k^*(\omega)$,~$x,y\in\cD_k$ and $\tau_E\PAR{x-y}\le \de_k$. \\

Let us recall that by Lemma 2.2
of \cite{OSSRF}, there exists $K_E\ge 1$ such that for all $x, y \in \mathbb{R}^d$
$$\tau_E(x+y)\le K_E\left(\tau_E(x)+\tau_E(y)\right).$$
Let us denote $F(h)=\sqrt{2dC_{\varphi}}\,h\ABS{\log h}^{(1+\eta)/\alpha+1/2}$ and choose $k_0\in\dN$ such that  $2^{k_0}\delta_{k_0+1}>3K_E^2$ and~$F$ is increasing on $(0,  \delta_{k_0}]$. Up to change $k^*\PAR{\omega}$, we can assume that $k^*\PAR{\omega}\ge k_0$. \\

Let $x,y\in\mathcal{D}$ with $x\neq y$ such that $3K_E^2\tau_E\PAR{x-y}\le \delta_{k^*(\omega)}$. Then there exists a unique $k\ge k^*(\omega)$ such that $\delta_{k+1}<3K_E^2\tau_E\PAR{x-y}\le\delta_k$. 
Since $x, y \in \mathcal{D}$, there exists $n\ge k+1$ such that  $x,y\in\cD_n$. Moreover, by Step 1, for $j=k,\ldots, n-1$, we can choose $x^{(j)},y^{(j)}\in\cD_j$ such that 
$$
\tau_E\left(x-x^{(j)}\right)\le 2^{-j} \ \textrm{and} \ \tau_E\left(y-y^{(j)}\right)\le 2^{-j}.
$$
By construction $\tau_E\left(x^{(k)}-y^{(k)}\right)\le K_E^2\left(\tau_E(x-y)+22^{-k}\right).$ 
Let us point out that since $k\ge k_0$, $2^k\delta_{k+1}\ge 2^{k_0}\delta_{k_0+1}>3K_E^2$. Therefore, one easily sees that 
$$
\tau_E\left(x^{(k)}-y^{(k)}\right)\le 3K_E^2\tau_E\PAR{x-y}. 
$$
Since $3K_E^2\tau_E(x-y)\le \delta_k$ we obtain by Step 3 on the one hand that
$$
\ABS{X_\al\PAR{x^{(k)}}-X_\al\PAR{y^{(k)}}}\le 
F\left(\tau_E\left(x^{(k)}-y^{(k)}\right)\right).$$
On the other hand we can write 
$$
X_\al(x)-X_\al\PAR{x^{(k)}}=\sum_{j=k}^{n-1}\PAR{X_\al\PAR{x^{(j+1)}}-X_\al\PAR{x^{(j)}}}
$$
with $\tau_E\left(x^{(j+1)}-x^{(j)}\right)\le 3K_E^22^{-(j+1)}\le \delta_{j+1}$ since $j\ge k_0$. 
Moreover, note that $x^{(j)}\in\cD_j\subset\cD_{j+1}$ and then  by Step 3
$$\ABS{X_\al(x)-X_\al\PAR{x^{(k)}}}\le \sum_{j=k}^{n-1}F\left(\tau_E\left(x^{(j+1)}-x^{(j)}\right)\right)\le CF(\delta_{k+1}),$$
where $C=\displaystyle\sum_{j=0}^{+\infty}\left(j+1\right)^{(1+\eta)/\alpha+1/2}\delta_j<+\infty$.
Using similar computations for $X_\al(y)-X_\al\PAR{y^{(k)}}$, we obtain that
\begin{eqnarray*}
\ABS{X_\al(x)-X_\al(y)}&\le& F\left(\tau_E\left(x^{(k)}-y^{(k)}\right)\right)+2CF\PAR{\delta_{k+1}}\\
&\le &(1+2C)F\left(3K_E^2\tau_E(x-y)\right).
\end{eqnarray*}
Then one can find a constant $C>0$ such that for $3K_E^2\tau_E\PAR{x-y}\le \delta_{k^*(\omega)}$
\begin{equation}\label{Majodx}
\ABS{X_\al(x)-X_\al(y)}\le C\tau_E(x-y)\ABS{\log \tau_E(x-y)}^{(1+\eta)/\alpha+1/2}.
\end{equation}
\\
We now give a modification of $X_\al$. For  $x\in \mathcal{D}$, we set
$$
X^*_\al\PAR{x}\PAR{\omega}=X_\al\PAR{x}\PAR{\omega}. 
$$
 For $x\in K$ let $x^{(n)}\in\cD$ such that
 $\lim_{n\to+\infty}x^{(n)}=x$. In view of~\eqref{Majodx}, $\PAR{X^*_\al\PAR{x^{(n)}}\PAR{\omega}}_n$ is  a Cauchy sequence  and then converges. We  set 
$$
X_\al^*\PAR{x}\PAR{\omega}=\lim_{n\to+\infty}X_\al^*\PAR{x^{(n)}}\PAR{\omega}.
$$
Remark that this limit does not depend on the choice of $\PAR{x^{(n)}}$. 
Moreover, since $X_\al$ is stochastically continuous, $X_\al^*$ is a modification of $X_\al$. 

By continuity of $\tau_E$, we easily see that 
$$
\ABS{X_\al^*\PAR{x}\PAR{\omega}-X_\al^*\PAR{x'}\PAR{\omega}}\le C\tau_E(x-y)\ABS{\log \tau_E(x-y)}^{(1+\eta)/\alpha+1/2}
$$
as soon as $3K_E^2\tau_E\PAR{x-y}< \delta_{k^*(\omega)},$
which concludes the proof. 
\end{proof}

 Following the same lines  as the proof of Theorem \ref{ModCont} we obtain
a similar result in the Gaussian case ($\al=2$) for more general fields. Let us remark that $Y_\alpha$ is not defined for $\alpha=2$. However, in Step 2 of the proof,  let us replace $Y_\alpha$ by 
$X$ a centered Gaussian random field and $v^2\PAR{\PAR{x,y}\left|\right.\PAR{T_n,\xi_n}_n}$ by the variance
of $X(x)-X(y)$
$$v^2\PAR{\PAR{x,y}}=\mathbb{E}\PAR{\PAR{X\PAR{x}-X\PAR{y}}^2}.$$
Furthermore let us replace Step 3 by the assumption that for some $\beta\in\mathbb{R}$ and $\delta>0$
there exists a finite constant $C>0$ such that for $x, y\in K$ with $\tau_E(x-y)\le \delta$
\begin{equation}\label{majogaussien}
\mathbb{E}\PAR{\PAR{X\PAR{x}-X\PAR{y}}^2}\le C\tau_E(x-y)^2\ABS{\log \tau_E(x-y)}^{\beta}.
\end{equation}
Then Step 1, Step 2 and Step 4 yields the following proposition.

\begin{prop}\label{ModContg} Let $X$ be a centered Gaussian random field satisfying \eqref{majogaussien}.
There exists  a modification $X^*$ of $X$ on $K$ such that 
\begin{equation}\label{moduleg}
\lim_{\delta\downarrow 0}\underset{\underset{\|x-y\|\le \delta}{x,y\in K}}{\sup}\frac{\ABS{X^*(x)-X^*(y)}}{\tau_E(x-y)\ABS{\log \tau_E(x-y)}^{1/2+\beta+\varepsilon}}=0  \mbox{ a.s.}
\end{equation}
for any $\varepsilon>0$. 
\end{prop}
Let us point out that if $X_2$ is an operator scaling Gaussian random field as defined in \cite{OSSRF}, then
$$\mathbb{E}\PAR{\PAR{X_2\PAR{x}-X_2\PAR{y}}^2}=\tau_E\PAR{x-y}^2\mathbb{E}\PAR{X_2\PAR{\ell_E(x-y)}^2},$$
and $X_2$ satisfies \eqref{majogaussien} with $\beta=0$ by $\PAR{5.2}$ of \cite{OSSRF}.
Therefore this result is more precise than one could expect from the Theorem \ref{ModCont},
replacing  $\al$ by $2$.

Let us also mention that Marianne Clausel gives a different proof of a similar result for some Gaussian operator scaling random fields with stationary increments in \cite{Marianne}.\\

Now, as in \cite{OSSRF}, we are looking for global and directional   H\"older critical exponents of the harmonizable stable  random field $X_\al$. These exponents have been introduced in \cite{Bonami} in the Gaussian realm but can be defined for any random field. Following Definition 5.1 of \cite{OSSRF}, $H\in (0,1)$ is
said to be the  H\"older
 critical exponent of a random field $\PAR{X(x)}_{x\in\rd}$ if there exists a modification $X^*$ of $X$ such that for any $s\in (0,H)$, the sample paths of $X^*$
satisfy almost surely a
 uniform H{\"o}lder condition of order $s$ on  $K$,
 that is
  there exists a finite positive random
 variable $A$ such that almost surely 
\begin{equation}\label{Holderu}
\left|X^*(x)-X^*(y)\right|\le A\|x-y\|^{s}\quad\text{for all } x,y \in
 K
 \end{equation}
while for any $s\in (H,1)$, almost surely \eqref{Holderu} fails. Note that  the  H\"older
 critical exponent, if exists, is well defined since any continuous modification of $X$ and $X^*$ are indistinguishable. Moreover, according to Definition 5.3 of
\cite{OSSRF} we say that  $X$ admits $H(u)$ as directional regularity in
direction $u\in S^{d-1}$, with $S^{d-1}$ the Euclidean unit sphere, if the process
$\PAR{X(tu)}_{t\in\mathbb{R}}$ admits  $H(u)$ as H\"older critical exponent on $K\cap \mathbb{R}u$.\\

For  all  $j=1,\ldots,m$ we set $K_j=K\cap \bigoplus_{k=1}^jW_k$. Let us remark that $K_m=K$.

\begin{cor}\label{ModContnorm} Let $\alpha\in (0,2]$. There exists  a modification $X^*_\al$ of $X_\al$  on $K$ such that 
for all  $j=1,\ldots,m$
$$
\lim_{\delta\downarrow 0}\underset{\underset{\|x-y\|\le \delta}{x,y\in K_j}}{\sup}\frac{\ABS{X^*_{\al}(x)-X^*_\al(y)}}{\|x-y\|^{H_j}\ABS{\log \NRM{x-y}}^{H_j(p_j-1)+\beta+1/2+\varepsilon}}=0  \mbox{ a.s.}
$$
for any $\varepsilon>0$, where $p_j=\underset{1\le k\le j}{\max}l_k$, $\beta=1/\alpha$ if $\alpha\neq 2$ and $\beta=0$ if $\alpha=2$. 
\end{cor}
\begin{proof}It follows from Theorem \ref{ModCont} and Corollary \ref{taulocal}, since $a_j \le a_d$ for any $j=1\ldots d$.
\end{proof}

\begin{cor}\label{Holder} Let $\alpha\in (0,2]$.
The random field $X^*_{\alpha}$  has locally $H$-H\"older sample paths on $\rd$ for every $H\in(0,H_m)$.
\end{cor}

Now let us give the  directional and global H\"older critical exponents of $X_{\alpha}$.

\begin{prop}\label{regularite directionnelle} The random field $X_{\alpha}$  admits $H_m$ as H\"older critical
exponent.\\ Moreover, for any $j=1,\ldots,m$, for any direction
$u\in W_j$, the field $X_\al$  admits $H_j$ as directional regularity in the direction
$u$.
\end{prop}
\begin{proof} For $Z$ a real S$\alpha$S random variable we let
$$\|Z\|_{\alpha}=\PAR{-\log\PAR{\mathbb{E}\PAR{\exp\PAR{iZ}}}}^{1/\alpha}.$$
Then, for any $x, y\in \rd$,
$$
\|X_\al^*(x)-X_\al^*(y)\|_{\alpha}=C_{\alpha}\PAR{\ell_E(x-y)}\tau_E(x-y)$$
where for all $\theta\in S_E$ $$C_{\alpha}\PAR{\theta}=\PAR{c_{\alpha}\int_{\rd}\left|\mbox{e}^{i\langle\theta,\xi\rangle}-1\right|^{\alpha}
\psi(\xi)^{-\alpha-q}d\xi}^{1/\alpha} \mbox{ and } c_{\alpha}=\frac{1}{2\pi}\int_0^{\pi}\ABS{\cos(t)}^\al dt.
$$
From Lebesgue's Theorem, the function $C_{\alpha}$ is continuous on the compact set $S_E$ with positive values. Let us denote 
$m_{\alpha}=\underset{\theta\in S_E}{\min}C_{\alpha}\PAR{\theta}>0$. According to~\eqref{control WjE} in  Corollary~\ref{tauglobal}, for any
$j=1,\ldots,m$ and $u$ a direction in $W_j$,
$$
\|X_\al^*(tu)-X_\al^*(su)\|_{\alpha}\ge m_{\alpha}c_1\left|t-s\right|^{H_j}\left|\log \left|t-s\right|\right|^{-(l_j-1)H_j}.
$$
Therefore, for any $s>H_j$, it implies that $\frac{X_\al^*(tu)-X_\al^*(su)}{\ABS{t-s}^s}$ is almost surely unbounded as $\ABS{t-s}\downarrow 0$ so~\eqref{Holderu} fails almost surely  on $K\cap \mathbb{R}u$.\\
Moreover, 
Corollary~\ref{ModContnorm} implies that $\PAR{X_\al^*(tu)}_{t\in\mathbb{R}}$ satisfies \eqref{Holderu} on $K\cap \mathbb{R}u$ for any $s<H_j$ and thus~$H_j$ is the directional regularity of $X_\al$ in the direction $u$.\\
Moreover, one can find a direction $u\in S^{d-1}$ in which almost surely $\PAR{X_\al^*(tu)}_{t\in\mathbb{R}}$ does not satisfy \eqref{Holderu} on $K\cap \mathbb{R}u$ for any $s>H_m$. Therefore, almost surely $\PAR{X_\al^*(x)}_{x\in\mathbb{R}^d}$
can not satisfy \eqref{Holderu} on $K$ for any~$s>H_m$. Then, by Corollary \ref{Holder} 
$X_\al$  admits $H_m$ as H\"older critical
exponent.\\
\end{proof}

\begin{rem}{\rm In the diagonalizable case (see Example \ref{diago}), the $W_j$ are the eigenspaces associated with the eigenvalues of $E$. In particular, for $\theta_j$ an eigenvector associated with the eigenvalue $\lambda_j=a_j$, the critical H\"older exponent in direction $\theta_j$ is $H_j=1/a_j$. }
\end{rem}

Proposition \ref{regularite directionnelle} compared to Theorem 5.4 of \cite{OSSRF} shows that the operator scaling stable field,
defined through an harmonizable representation share the same sample paths properties as the Gaussian ones. Therefore it is natural to have also the same result of  Theorem 5.6  \cite{OSSRF} for the  box- and the
Hausdorff-dimensions of their graphs on a compact set. We also refer to
Falconer \cite{Falconer} for the definitions and properties of
box- and the Hausdorff-dimension and keep the notations of \cite{OSSRF}.
More precisely, we fix a compact set
$K\subset\rd$ and consider
$\mathcal G(X_\al^*)(\omega)=\{(x,X_\al^*(x)(\omega)); x\in K\}$ the graph of
a realization of the field $X_\al^*$ over the compact $K$. We  denote
$\mbox{dim}_{\mathcal H}\mathcal G(X_\al^*)$, resp $\mbox{dim}_{\mathcal
B}\mathcal G(X_\al^*)$, the Hausdorff-dimension and the box-dimension of
$\mathcal G(X^*_\al)$, respectively.
\begin{prop}\label{XHD} Almost
surely
$$\dim_{\mathcal H}\mathcal G(X^*_\al)=\dim_{\mathcal B}\mathcal G(X^*_\al) = d+1-H_m.$$
\end{prop}

\begin{proof}The proof is very similar to those of Theorem 5.6  \cite{OSSRF}. It also uses same kinds of arguments as in
\cite{Benassi}. For sake of completeness we recall the main ideas. Corollary \ref{Holder} allows as usual
to state the upper bound
$$\mbox{ dim}_{\mathcal H}\mathcal G(X_\al^*)\leq\overline{\mbox{ dim}_{\mathcal B}}\mathcal G(X_\al^*) \le d+1-H_m,\text{ a.s.}$$
where $\overline{\mbox{ dim}_{\mathcal B}}$ denotes the upper
box-dimension. The lower bound will also follows from Frostman criterion (Theorem 4.13 (a) in \cite{Falconer}). One has to prove that the
integral $I_s$
$$I_s=\int_{K\times K}\mathbb{E}\left[\left((X^*_\al(x)-X^*_\al(y))^2+\| x-y\|^2\right)^{-s/2}\right]\,dx\,dy,$$
is finite to get that  almost surely $\mbox{ dim}_{\mathcal H}\mathcal
G(X_\al^*)\ge s$. In our case, the fundamental lemma of \cite{AyacheR} allows us to write this integral using the characteristic function of the S$\alpha$S field $X_\al^*$. Actually, when one remarks that, using Fourier-inversion,
$(\xi^2+1)^{-s/2}=\frac{1}{2\pi}\int_{\mathbb{R}}e^{i\xi t}f_s(t)dt$, where $f_s\in L^{\infty}(\mathbb{R})\cap L^1(\mathbb{R})$,
 one gets
$$I_s=\int_{K\times K}\left(\frac{1}{2\pi}\| x-y\|^{-s}\int_{\mathbb{R}}
e^{-|t|^\alpha\frac{\left\|X_\al^*(x)-X_\al^*(y)\right\|_{\alpha}^{\alpha}}{\|
x-y\|^\alpha}}f_s(t)dt \right)dx \,dy.$$ 
By a change of variables, as $f_s\in L^{\infty}(\mathbb{R})$, one can find $C>0$ such that
\begin{eqnarray*}
I_s&\le & C\int_{K\times K}\| x-y\|^{1-s}\left\|X_\al^*(x)-X_\al^*(y)\right\|_{\alpha}^{-1}dx\, dy\\
&\le & Cm_{\alpha}^{-1}\int_{K\times K}\| x-y\|^{1-s}\tau_E(x-y)^{-1}dx\, dy,
\end{eqnarray*}
where $m_{\alpha}=\underset{\theta\in S_E}{\min}\PAR{c_{\alpha}\int_{\rd}\left|\mbox{e}^{i\langle\theta,\xi\rangle}-1\right|^{\alpha}
\psi(\xi)^{-\alpha-q}d\xi}^{1/\alpha}$ as introduced in the proof of Proposition \ref{regularite directionnelle}.
Since $\int_{K\times K}\| x-y\|^{1-s}\tau_E(x-y)^{-1}dx \,dy$ is proved to be finite as soon as $s<d+1-H_m$ in   \cite{OSSRF},
 $$\underline{\mbox{ dim}_{\mathcal B}}\mathcal
G(X^*_\al)\ge \mbox{ dim}_{\mathcal H}\mathcal
G(X^*_\al)\ge d+1-H_m \mbox{ a.s. }$$  and the proof is complete.\\
\end{proof}

Harmonizable operator scaling stable random fields share many properties with Gaussian operator random fields. In particular, they have locally H\"older sample paths and  critical directional H\"older exponent   depending  on the directions. In the Gaussian case ($\al=2$),~\cite{OSSRF} establishes such properties in the framework of  harmonizable and moving average Gaussian operator scaling random field. 
However, for stable laws, harmonizable and moving average representations do not have the same behavior as we see  in the next section. 
\section{Moving average representation}
\label{MASE}

 Let us recall the definition  of moving average operator scaling stable random fields introduced in~\cite{OSSRF}. Let $0<\alpha\le  2$. We consider $M_\alpha(dy)$ an independently scattered $S\alpha S$ random measure on $\rd$ with
Lebesgue control measure, see~\cite{sam:taqqu} for details on such random measures. Let us recall that, as before,  $q=\mbox{trace}(E)$. Let $\varphi : \rd\to[0,\infty)$ be a continuous $E$-homogeneous function. We assume moreover that there exists $\beta>1$ such that $\varphi$ is {\it $\left(\beta, E\right)$-admissible}. According to Definition 2.7 of \cite{OSSRF} it means that $\varphi(x)\neq 0$ for $x\neq 0$ and that  for any $0<A<B$ there exists
a  constant $C>0$ such that, for $A\le \|y\|\le B$,
\begin{equation}\label{admissible}
 |\varphi(x+y)-\varphi(y)|\leq C\tau_E(x)^{\beta}
 \end{equation}
 holds for any $\tau_E(x)\le 1$.\\

 \begin{defi}
Since $\beta>1$, the $\alpha$-stable random field
\begin{equation}\label{moving average}
Z_\al(x)=\int_\rd\Bigl(\varphi(x-y)^{1-q/\alpha}-\varphi(-y)^{1-q/\alpha}\Bigr)\,M_\alpha(dy)\ ,x\in\rd.
\end{equation}
is well defined and called moving average operator scaling stable random field. 
\end{defi}
From Theorem 3.1 and  Corollary 3.2 of \cite{OSSRF},  it is stochastically continuous, has stationary increments and satisfies the following operator scaling property 
$$\forall c>0, \ \BRA{Z_\al(c^Ex) ; x\in\rd}\dlim\BRA{cZ_\al(x) ; x\in\rd},$$ as  the harmonizable field $X_\al$. \\

In the Gaussian case ($\al=2$), the variogramme of $Z_2$ is similar to the one of the harmonizable field $X_2$. Then,~\cite{OSSRF} proves that $Z_2$ and $X_2$ admit the same critical H\"older sample paths properties. However, when $\al\in(0,2)$, let us recall that  moving average self-similar $\al$-stable random motions  does not have in general continuous sample paths (see \cite{sam:taqqu}). The next proposition states the same property for $Z_\al$.

\begin{prop} Assume $\al\in(0,2)$ and $d\ge 2$. Then, the random field $Z_\al$ is almost surely unbounded on every open ball.
\end{prop}

\begin{proof} 
Let us remark that $\varphi(0)=0$ by continuity and $E$-homogeneity and $q=\underset{j=1}{\overset{d}{\sum}}a_j>d>\alpha$, as soon as $d\ge 2$. Then, for any $U$ open set, since $U^*=U\cap \mathbb{Q}^d$ is dense in $U$, for any $y\in U$
$$
f^*\PAR{U^*,y}=\underset{x\in U^*}{\sup}\ABS{\varphi(x-y)^{1-q/\alpha}-\varphi(-y)^{1-q/\alpha}}=+\infty.
$$
Then $\int_{\rd}f^*\PAR{U^*,y}^{\alpha}dy=+\infty$ and the necessary condition for sample boundedness (10.2.14)
of Theorem 10.2.3 p.450 of \cite{sam:taqqu} fails. We conclude the proof by Corollary 10.2.4 of \cite{sam:taqqu}.
\end{proof}

\section{Appendix}
\begin{proof}[Proof of lemma~\ref{Jordan power}]

The lower bound is straightforward. Actually, for any $t>0$, $t^{\lambda}$ is an eigenvalue of the matrix $t^J$  and therefore $t^{a}=\left|t^{\lambda}\right|\le \|t^{J}\|$. 

Let us prove the upper bound. First, let us assume that 
$J$ is a Jordan cell matrix of size $l$. In this case $\lambda=a\in\dR$ and 
$$t^{J}= t^{a}\left(\begin{array}{cccc} 1& 0 &\ldots & 0\\
\log t & 1 & 0&\vdots\\
\vdots&\ddots& \ddots& 0\\
\frac{(\log t)^{l-1}}{(l-1)!}&\ldots& \log t & 1\end{array}\right).$$
Let us recall that the norm defined 
for a matrix $A=(a_{ij})_{1\le i,j,\le d}$ by $\|A\|_{\infty}=\underset{1\le i\le d}{\max}\sum_{j=1}^d|a_{ij}|$ is the subordinated norm of the infinite norm $\|x\|_{\infty}=\underset{1\le i\le d}{\max}|x_i|$ for $x\in \mathbb{R}^d$. By definition, we   can deduce that $\|t^{J}\|_{\infty}=t^{a}\underset{j=0}{\overset{l-1}{\sum}}\frac{\left|\log t\right|^{j}}{j!}$.
Then,
$$
\NRM{t^J}\le\sqrt{l}\NRM{t^J}_\infty\le \sqrt{l} t^{a}\ABS{\log\PAR{t}}^{l-1}\underset{j=0}{\overset{l-1}{\sum}}\frac{1}{j!}  
$$
for any $t\in (0,e^{-1}]\cup[e,+\infty)$. Therefore,  for any $t\in (0,e^{-1}]\cup[e,+\infty)$, we have
$$ \|t^{J}\|\le \sqrt{l}\,e\, t^{a}\left|\log t\right|^{l-1}.$$

 In the second case, let us assume that $J$ is a block of the form~\eqref{Jordan block} of size $2l$ associated with the 
eigenvalue $\lambda=a+ib$ for $b\neq 0$. Then $t^{J}=t^{a}R(t)N(t)$ where
$$R(t)=
\left(\begin{array}{cccc} R_{b}(t)& 0 &\ldots & 0\\
0 & R_{b}(t) & 0&\vdots\\
\vdots&\ddots& \ddots& 0\\
0&\ldots& 0 & R_{b}(t)\end{array}\right)
\mbox{ with } R_{b}(t)=\left(\begin{array}{cc} \cos(b\log t) &-\sin(b\log t)\\
\sin(b\log t)&  \cos(b\log t)\end{array}\right),$$
and
$$N(t)=
\left(\begin{array}{cccc} I_2& 0 &\ldots & 0\\
N_1(t) & I_2 & 0&\vdots\\
\vdots&\ddots& \ddots& 0\\
N_{l-1}(t)&\ldots& N_1(t) & I_2\end{array}\right)\mbox{ with } N_j(t)=\left(\begin{array}{cc} \frac{\left|\log t\right|^{j}}{j!} &0\\
0&  \frac{\left|\log t\right|^{j}}{j!}\end{array}\right).$$
Hence,
$$
\NRM{t^J}\le t^\al\NRM{R\PAR{t}}\NRM{N\PAR{t}}.
$$
Since $R\PAR{t}$ is an orthogonal matrix, $\NRM{R\PAR{t}}=1$. Furthermore, $N\PAR{t}$ is a $(2l)\times (2l)$ matrix and 
$$
\NRM{N\PAR{t}}\le\sqrt{2l}\NRM{N(t)}_\infty=\sqrt{2l} \underset{j=0}{\overset{l-1}{\sum}}\frac{\left|\log t\right|^{j}}{j!}
$$
Therefore, 
we also obtain that  
$$ \|t^{J}\|\le  \sqrt{2l} \, e t^{a}\left|\log t\right|^{l-1}$$ for any $t\in (0,e^{-1}]\cup[e,+\infty)$.\\
\end{proof}
\begin{proof}[Proof of Proposition~\ref{tauglobal}]
Let $r\in (0,1)$. 
One can find $r_E\in(0,r)$ such that for any $\|x\|\le r_E$ we have $\tau_E(x)\le e^{-1}$. Let $x\in \bigoplus_{k=j_0}^j W_k$ with $\|x\|\le r_E$.  Then $x=\tau_E(x)^E\ell_E(x)$ and $l_E\PAR{x}\in \bigoplus_{k=j_0}^j W_k$. \\

We first establish the lower bound of Proposition~\ref{tauglobal}. Let us write $l_E\PAR{x}=\sum_{k=1}^m l_k\PAR{x}$ where each $l_k\PAR{x}\in W_k$. 
   Let $L_k$ be the coordinates of $l_k\PAR{x}$ in the basis $\PAR{f_{\sum_{i=1}^{k-1}\tilde{l_i}+1},\ldots,f_{\sum_{i=1}^{k}\tilde{l_i}}}$ of $W_k$.  Hence, by definition of $P$, 
$$
P^{-1}l_E\PAR{x}=\PAR{\begin{array}{c}L_{1}\\\vdots\\L_m\end{array}}
 \ 
and 
 \ 
x=\tau_E\PAR{x}^El_E\PAR{x}=P\PAR{\begin{array}{c}\tau_E\PAR{x}^{J_1}L_1\\\vdots\\\tau_E\PAR{x}^{J_m}L_m\end{array}}.
$$
Since $l_E\PAR{x}\in \bigoplus_{k=j_0}^j W_k$, $L_k= 0$ for $k\notin\{j_0,\ldots,j\}$,
$$
\NRM{x}\le \NRM{P} \PAR{\displaystyle\sum_{k=j_0}^j\NRM{\tau_E\PAR{x}^{J_k}L_k}^2}^{1/2}\\\\
\le \NRM{P}\PAR{\sum_{k=j_0}^j\NRM{\tau_E\PAR{x}^{J_k}}^2\NRM{L_k}^2}^{1/2}\\
$$
By Lemma~\ref{Jordan power}, 
$$\NRM{x}
\le\sqrt{2}e\NRM{P}\PAR{\sum_{k=j_0}^j l_k \tau_E\PAR{x}^{2a_k}\ABS{\log\tau_E\PAR{x}}^{2\PAR{l_k-1}}\NRM{L_k}^2}^{1/2}
$$
since $\tau_E\PAR{x}\le 1/e$. Hence, since $\tau_E\PAR{x}\le 1$, $\ABS{\log\tau_E\PAR{x}}\ge 1$, $a_k\ge a_{j_0}$ and $l_k\le p_{j_0,j}=\displaystyle\max_{j_0\le i\le j}l_i\le d$,
$$\begin{array}{rcl}
\NRM{x}
&\le&\sqrt{2d}e\NRM{P}\tau_E\PAR{x}^{a_{j_0}}\ABS{\log\tau_E\PAR{x}}^{\PAR{p_{j_0,j}-1}}\PAR{\sum_{k=j_0}^j \NRM{L_k}^2}^{1/2}\\&\le&\sqrt{2d}e\NRM{P}\tau_E\PAR{x}^{a_{j_0}}\ABS{\log\tau_E\PAR{x}}^{\PAR{p_{j_0,j}-1}}\NRM{P^{-1}l_E\PAR{x}}\\
\end{array}
$$
Then, 
\begin{equation}
\label{mino1}
\NRM{x}\le \sqrt{2d}eM_E\NRM{P}\NRM{P^{-1}}\tau_E\PAR{x}^{a_{j_0}}\ABS{\log\tau_E\PAR{x}}^{\PAR{p_{j_0,j}-1}}
\end{equation}
where $M_E$ is defined by~\eqref{ME}.
 
Let us take the logarithm of this equation. Choosing  $r_0<\min(1,r_E)$ small enough, one can find $C>0$ such that 
\begin{equation}\label{Log}
|\log \tau_E(x)|\le C |\log \|x\|| \mbox{ for } \|x\|<r_0. 
\end{equation}

Using in~\eqref{mino1}, we obtain that there exists $C>0$ such that for $\NRM{x}\le r_0$
$$
\NRM{x}^{H_{j_0}}\ABS{\log\NRM{x}}^{-H_{j_0}\PAR{ p_{j_0,j}-1}}\le C \tau_E\PAR{x}.
$$
 By continuity of the map, $$x\mapsto \NRM{x}^{H_{j_0}}\ABS{\log\NRM{x}}^{-H_{j_0}\PAR{ p_{j_0,j}-1}}\tau_E\PAR{x}^{-1}$$ on $\BRA{x\in\rd\,/\, 0<\NRM{x}<1}$, up to change $C$, the previous inequality holds for $\NRM{x}\le r$, which gives the lower bound in~\eqref{control uE}. \\

Let us now establish the upper bound in~\eqref{control uE}. We write $x=\sum_{k=1}^m x_k$ with each $x_k\in W_k$. We then denote by $X_k$ the coordinates  of $x_k$ in the basis $\PAR{f_{\sum_{i=1}^{k-1}\tilde{l_i}+1},\ldots,f_{\sum_{i=1}^{k}\tilde{l_i}}}$ of $W_k$. Since $x\in \bigoplus_{k=j_0}^jW_k$,
$
X_k=0$ for every $k\notin\{j_0,\ldots,j\}$. 
Hence, by definition of $P$, 
$$
P^{-1}x=\PAR{\begin{array}{c}X_1\\\vdots\\X_m\end{array}}
 \ 
\mbox{ and } 
 \ 
l_E\PAR{x}=\tau_E\PAR{x}^{-E}x=P\PAR{\begin{array}{c}\tau_E\PAR{x}^{-J_1}X_1\\\vdots\\\tau_E\PAR{x}^{-J_m}X_m\end{array}}.
$$ 
Then, $\displaystyle\NRM{l_E\PAR{x}}\le \NRM{P}\PAR{\sum_{k={j_0}}^j\NRM{\tau_E\PAR{x}^{-J_k}X_k}^2}^{1/2}$.\\

The Lemma~\ref{Jordan power} yields 
$$
\NRM{l_E\PAR{x}}\le \sqrt{2}e\NRM{P}\PAR{\sum_{j=j_0}^j l_k \tau_E\PAR{x}^{-2a_k}\ABS{\log\tau_E\PAR{x}}^{2\PAR{l_k-1}}\NRM{X_k}^2}^{1/2}
$$
since $\tau_E\PAR{x}^{-1}\ge e$. Hence, since $\tau_E\PAR{x}^{-1}\ge e>1$, $a_k\le a_j$ and $l_k\le p_{j_0,j}$, 
$$
\begin{array}{rcl}
0<m_E&\le& \sqrt{2d}e\NRM{P}\NRM{P^{-1}}\tau_E\PAR{x}^{-a_j}\ABS{\log\tau_E\PAR{x}}^{ p_{j_0,j}-1}\sqrt{\sum_{k=j_0}^j \NRM{X_k}^2}\\&\le &\sqrt{2d}e\NRM{P}\NRM{P^{-1}}\tau_E\PAR{x}^{-a_j}\ABS{\log\tau_E\PAR{x}}^{ p_{j_0,j}-1}\NRM{P^{-1}x}.
\end{array}
$$
Then, using~\eqref{Log} and $\NRM{P^{-1}x}\le\NRM{P^{-1}}\NRM{x}$, there exists a constant $C>0$ such that 
$$
\tau_E\PAR{x}<C\NRM{x}^{H_j}\ABS{\log\NRM{x}}^{H_j\PAR{ p_{j_0,j}-1}}
$$
for $\NRM{x}\le r_0$. By continuity of the map
$$
x\mapsto \frac{\tau_E\PAR{x}}{\NRM{x}^{H_j}\ABS{\log\NRM{x}}^{H_j\PAR{ p_{j_0,j}-1}}}
$$
on $\BRA{x\in\rd\,/\, 0<\NRM{x}<1}$, up to change $C$, the previous inequality holds for $\NRM{x}\le r$, which gives the upper bound in~\eqref{control uE} and concludes the proof. \\
\end{proof}

\begin{proof}[Proof of Lemma~\ref{esperance}] It is sufficient to consider
 $$
 I(h)=\dE\PAR{m\PAR{\xi_n}^{-2/\al}\min\PAR{M_E\left\|h^{E^t}\xi_n\right\|,2}^2\psi\PAR{\xi_n}^{-2-2q/\al}} 
 $$
By definition,
$$
\begin{array}{rcl}
\displaystyle I(h)
&=&\displaystyle\int_{\dR^d}m\PAR{\xi}^{1-2/\al}\psi\PAR{\xi}^{-2-2q/\al}\min\PAR{M_E\left\|h^{E^t}\xi\right\|,2}^2d\xi.
\end{array}
$$
Using the formula of integration in {\it polar coordinates} with respect to $E^t$, see Proposition \ref{CDV}, 
$$
I(h)=\int_{S_{E^t}}\int_{0}^{+\infty}m\PAR{r^{E^t}\theta}^{1-2/\al}\psi\PAR{r^{E^t}\theta}^{-2-2q/\al}\min\PAR{M_E\left\|(hr)^{E^t}\theta\right\|,2}^2r^{q-1}dr\sigma_{E^t}\PAR{d\theta}.
$$
Since $\psi$ is $E^t$-homogeneous,
$$
\begin{array}{rcl}
I(h)&=&\displaystyle\int_{S_{E^t}}\int_{0}^{+\infty}m\PAR{r^{E^t}\theta}^{1-2/\al}\psi\PAR{\theta}^{-2-2q/\al}\min\PAR{M_E\left\|(hr)^{E^t}\theta\right\|,2}^2r^{-2+q-1-2q/\al}dr\sigma_{E^t}\PAR{d\theta}\\\\
&=&c_{\eta}^{1-2/\alpha}\displaystyle\int_{S_{E^t}}\int_{0}^{+\infty}\psi\PAR{\theta}^{-2-2q/\al}\min\PAR{M_E\left\|(hr)^{E^t}\theta\right\|,2}^2r^{-3}\ABS{\log(r)}^{\PAR{1+\eta}\PAR{2/\al-1}}dr\sigma_{E^t}\PAR{d\theta}.
\end{array}
$$
By the change of variable $\rho=hr$, $I(h)$ is equal to
$$
c_{\eta}^{1-2/\alpha}h^{2}\int_{S_{E^t}}\int_{0}^{+\infty}\psi\PAR{\theta}^{-2-2q/\al}\min\PAR{M_E\left\|\rho^{E^t}\theta\right\|,2}^2\rho^{-3}\ABS{\log\PAR{\frac{\rho}{h}}}^{\PAR{1+\eta}\PAR{2/\al-1}}dr\sigma_{E^t}\PAR{d\theta}.
$$
For any $\ga\in(0,1)$, there exists $A_\ga$ such that for every $\rho>0$ and every  $h\le 1-\ga$, 
$$
\ABS{\log\PAR{\frac{\rho}{h}}}=\ABS{\log\PAR{\rho}-\log\PAR{h}}\le A_\ga\ABS{\ABS{\log(\rho)}+1}\ABS{\log\PAR{h}}. 
$$
Since $2/\al>1$, 
$$
I(h)\le A_\ga^{2/\al-1}c_{\eta}^{1-2/\alpha}h^{2}\ABS{\log\PAR{h}}^{\PAR{1+\eta}\PAR{2/\al-1}}\PAR{I_1+I_2}
$$
with 
$$
I_1=4\int_{S_{E^t}}\psi\PAR{\theta}^{-2-2q/\al}\sigma_{E^t}\PAR{d\theta}\int_{1}^{+\infty}\rho^{-3}\ABS{\ABS{\log\PAR{\rho}}+1}^{\PAR{1+\eta}\PAR{2/ \al-1}}d\rho
$$
and 
$$
I_2=M_E^{2}M_{E^t}^{2}\int_{S_{E^t}}\psi\PAR{\theta}^{-2-2q/\al}\sigma_{E^t}\PAR{d\theta}\int_{0}^{1}\NRM{\rho^{E^t}}^{2}\rho^{-3}\ABS{\ABS{\log\PAR{\rho}}+1}^{\PAR{1+\eta}\PAR{2/ \al-1}}d\rho,
$$
where $M_E$ and $M_{E^t}$ are defined by~\eqref{ME}.
Since $\psi$ is continuous with positive value on the compact set~$S_{E^t}$, 
$$
\int_{S_{E^t}}\psi\PAR{\theta}^{-2-2q/\al}\sigma\PAR{d\theta}<+\infty.
$$
Hence $I_1<+\infty$.\\
It follows from Proposition~\ref{tauglobal} that for any $\de'\in(0,1)$, there exists a constant $c_\de'>0$ such that  
$$\|\rho^{E^t}\|\leq
C\rho^{a_1}\ABS{\log\ABS{\rho}}^{l-1}
$$ for all $\rho\leq \de'$. Hence, since $a_1>1$,  
$$
\int_{0}^{1}\NRM{\rho^{E^t}}^{2}\rho^{-3}\ABS{\ABS{\log\PAR{\rho}}+1}^{\PAR{1+\eta}\PAR{2/ \al-1}}d\rho<+\infty
$$
and $I_2<+\infty$, which concludes the proof.
\end{proof}

\bibliographystyle{plain}
\bibliography{Biblio}

\end{document}